\documentstyle{amsppt}
\magnification\magstep1
\baselineskip=15pt
\input btxmac.tex               
\bibliographystyle{plain}
\input epsf
\input diagrams
\NoRunningHeads
\hsize6.9truein
\pageheight{23 truecm}
\def\dim{\text{dim}}

\def\Hom{\text{Hom}}

\def\Sym{\text{Sym}}

\def\fb#1.{\underset #1 \to \times}

\def\Hom{\operatorname{Hom}}

\def\dim{\operatorname{dim}}
\def\Pic{\operatorname{Pic}}

\def\map#1.#2.{#1 \longrightarrow #2}
\def\ses#1.#2.#3.{0\longrightarrow #1 \longrightarrow #2 \longrightarrow #3 
\longrightarrow 0}
\def\es#1.#2.#3.{#1 \longrightarrow #2 \longrightarrow #3}
\def\ring#1.{\Cal O_{#1}}
\def\rmap#1.#2.{#1 \dasharrow #2}
\def\pr #1.{\Bbb P^{#1}}
\def\pl #1.#2.{#1^{\otimes #2}}
\def\proj#1.{\Bbb P(#1)}

\def\uloopr#1{\ar@'{@+{[0,0]+(-4,5)} @+{[0,0]+(0,10)}
@+{[0,0]+(4,5)}}
  ^{#1}}
\def\dloopr#1{\ar@'{@+{[0,0]+(-4,-5)} @+{[0,0]+(0,-10)}
@+{[0,0]+(4,-5)}}
  _{#1}}
\def\rloopd#1{\ar@'{@+{[0,0]+(5,4)} @+{[0,0]+(10,0)}
@+{[0,0]+(5,-4)}}
  ^{#1}}
\def\lloopd#1{\ar@'{@+{[0,0]+(-5,4)} @+{[0,0]+(-10,0)}
@+{[0,0]+(-5,-4)}}
  _{#1}}

\long\def\ignore#1{}
\long\def\ignore#1{#1}

\TagsOnRight
\NoBlackBoxes

\topmatter
\title Birational unboundedness of $\Bbb{Q}$-Fano threefolds \endtitle
\author JIAYUAN LIN \endauthor
\abstract We prove that the family of $\Bbb{Q}$-Fano threefolds with  Picard number 
one is birationally unbounded.
\endabstract

\endtopmatter

\subhead { 1. Introduction }\endsubhead

\smallpagebreak

We work over an algebraically closed field of characteristic zero.  

\proclaim{Definition 1.1} Let $X$ be a projective variety,  $X$ is said to be a $\Bbb Q$-Fano 
variety, if
\roster
\item $X$ has $\Bbb Q$-factorial log terminal singularities, and
\item $-K_X$ is ample.
\endroster
\endproclaim

The reader is referred to [14] for the standard definitions of higher-dimensional geometry,
such as log terminal.  The main goal of this paper is to prove the following main theorem.

\proclaim{Theorem 1.2} The family of $\Bbb{Q}$-Fano threefolds with Picard number one is  
birationally unbounded.
\endproclaim

For the reader's convenience we state the definition of birationally bounded which is a
modification of the definition of bounded given by V. Alexeev in [1].

\proclaim{Definition 1.3} A class of varieties $\Cal B$ is birationally bounded if there
exists a morphism $f:\map \Cal X. \Cal S.$ between two varieties such that
every variety in $\Cal B$ is birational to one of the geometric fibres of $f$.
\endproclaim

Note that we do not require that every geometric fibre of $f$ is birational to a variety
in $\Cal B$, nor do we require that every geometric fibre of $f$ is birational to a unique
variety in $\Cal B$.

 $\Bbb Q$-Fano varieties play an important role in modern birational algebraic geometry.
In the Minimal Model Program (MMP) we have to allow varieties with certain mild
singularities.  $\Bbb Q$-Fano varieties appear naturally as one of the final results of
running the log MMP.  The following is an interesting and motivating conjecture proposed
independently by A. Borisov [2] and V. Alexeev [1].

\proclaim{Conjecture 1.4 (A.Borisov-V.Alexeev)} Fix $\epsilon >0$.  Then the family of 
all $\Bbb Q$-Fano varieties of a given dimension with log discrepancy greater than
$\epsilon$ is bounded.
\endproclaim

\remark {Remark} One cannot remove $\epsilon$ from the hypothesis.  Take, for example, 
the cone over a rational curve of degree $d$.  For every $d$ the corresponding cone has
log discrepancy $\frac{2}{d}$.  These cones are all $\Bbb Q$-Fano surfaces (also known as
log Del Pezzo surfaces) but they form an unbounded family.
\endremark

 As well as being an interesting conjecture in its own right, Conjecture 1.4 and many similar
conjectures play a pivotal role in an inductive approach to higher dimensional geometry.
Conjecture 1.4 also has important applications to the Sarkisov program.

 For these reasons there has been a considerable amount of work on this conjecture.
A.Nadel [16] and F. Campana [5] proved the boundedness of smooth Fano varieties with
Picard number one. J. Koll\'ar, Y. Miyaoka and S. Mori [11] proved the boundedness of
smooth Fano varieties with arbitrary Picard number in every dimension.  The case of $\Bbb
Q$-Fano threefolds with Picard number one and terminal singularities is due to Y. Kawamata
[9].  Z. Ran and H. Clemens [18] proved a boundedness theorem for Fano unipolar
$n$-dimensional varieties.  V. Alexeev [1] and V.V. Nikulin [17] proved the above
conjecture in dimension two.  A. Borisov and L. Borisov [4] gave a proof of Conjecture 1.4 in the
toric case.  A. Borisov [2], [3] proved that $\Bbb Q$-Fano threefolds of given index are
bounded. J. Koll\'ar, Y. Miyaoka, S. Mori and H. Takagi [12] proved that all $\Bbb Q$-Fano
threefolds with canonical singularities are bounded.  Recently J. M$^c$Kernan [15] proved
boundedness for log terminal Fano pairs of bounded index. Despite this, Conjecture 1.4 is
unresolved, even in dimension three.

 Consider cones over rational curves of degree $d$ as above.  Even though they form an
unbounded family, all of them are birational to $\pr 2.$.  More generally any quotient of
$\pr 2.$ by a finite group $G$ is a log Del Pezzo surface of Picard number one.  The
classification of log Del Pezzo surfaces has attracted considerable interest.  V. Alexeev
and V.V. Nikulin classified log Del Pezzo surfaces under certain conditions on the indices
of the singularities.  One of the key steps to prove that the family of smooth Fano
varieties is bounded is to bound the top self-intersection of $-K_X$.  Once one allows
singularities life is not so simple.  In fact S. Keel and J. M$^{c}$Kernan [10] proved
that the set
$$
\{\, K_S^2 \,|\, \text {$S$ is a toric surface of Picard number one} \,\}
$$
is dense in $\Bbb R_+$.  They classified all but a bounded family of
Picard number one log Del Pezzo surfaces.  Although there are very many log Del Pezzo
surfaces up to isomorphism, all of them are rational.  Thus the family of log Del Pezzo
surfaces is birationally bounded and in fact birationally bounded is considerably weaker
than bounded.

 It has been speculated that the same is true in the higher dimensional case, that is to
say it has been speculated that $n$-dimensional $\Bbb Q$-Fano varieties of Picard number
one are birationally bounded.  We find this is not the case.  As a consequence, one cannot
remove $\epsilon$ from the hypothesis of Conjecture 1.4, even if we replace bounded by birationally
bounded.

 It is also interesting to note that we cannot drop the condition that the singularities
are log terminal even in the case of surfaces.  Consider the family of surfaces with
ample anticanonical divisor and Picard number one.  This family is birationally unbounded.
For example take a cone over a curve of genus $g$ and degree greater than $2g-2$.

 The idea of the proof of Theorem 1.2 is to look at conic bundles contained in $\pr
2.$-bundles over $\pr 2.$, obtained as the blow ups of cones over Veronese embeddings of
$\pr 2.$.  Varying the degree of the embedding gives an infinite family and, using deep
results of Sarkisov concerning the classification of conic bundles, we prove that the
conic bundles so constructed are not birational to each other.  Finally we use the MMP and
some results of Iskovskikh [8], suitably modified, to prove that these conic bundles are
birationally unbounded.

\subhead Acknowledgements\endsubhead 

It is a great pleasure to express my appreciation to my advisor James M$^{c}$Kernan for
proposing this problem and for his help and encouragement. I am grateful to Xi Chen and
Mikhail Grinenko for their kind help.  I would also like to thank the referee for numerous
suggestions which have considerably improved the exposition.  This paper uses Paul
Taylor's commutative diagrams package written in \TeX.

In order to prove Theorem 1.2, we need some preliminary results.

\subhead {2. Some examples of rigid conic bundles}\endsubhead

\definition{Definition 2.1} A projective variety $X$ is said to be a conic bundle, if there
is a morphism $\pi:\map X.S.$, where $S$ is an irreducible projective variety and the
geometric generic fibre is a smooth rational curve.

  If there exists a rank $3$ locally free sheaf $\Cal E$ on $S$ such that $\pi$ factors
into a closed embedding $X \hookrightarrow \pr. (\Cal E) $ as the zero locus of a global
section $v$ of $ \Sym^{2} \Cal E$ and the natural projection $\map {\pr.(\Cal E)}. S.$,
then $X$ is a classic conic bundle.  The zero locus $\Delta(v)$ of $v$ is called
the discriminant locus.
\enddefinition

 Note that the discriminant locus $\Delta(v)$ of a classic conic bundle is a divisor and
that as a set it parametrises the singular fibres.  We adopt these properties of the
discriminant locus as the definition of the discriminant locus for an aribtrary conic
bundle.

\definition{Definition 2.2} Let $\pi:\map X.S.$ be a conic bundle. Then the 
discriminant locus $\Delta(\pi)$ is defined to be the codimension one part of $\pi(R)$,
where $R$ is the locus in $X$ where $\pi$ is not smooth.
\enddefinition

\example{Example} Let $S$ be any smooth surface and let $X$ be the ordinary blow up of 
$S\times \pr 1.$ along any point.  Then $\map X.S.$ is a conic bundle.  On the other hand
there is only one singular fibre.  Thus even though the singular locus is non-empty, the
discriminant locus is empty by our definition.
\endexample

\definition{Definition 2.3} A conic bundle $\pi:\map X.S.$ is called standard if
\roster
\item $X$ and $S$ are smooth projective varieties,
\item $\pi$ is a flat morphism, 
\item $\rho(X)=\rho(S)+1$, and
\item $-K_{X}$ is relative ample.
\endroster
\enddefinition

\definition{Definition 2.4} Assume $X$ has only $\Bbb Q$-factorial terminal singularities. 
A Mori fibre space is an extremal contraction $f:\map X.S.$ of fibre type. In other words
$f_* \ring X.= \ring S.$ and
\roster
\item $-K_X$ is relatively ample for $f$,
\item $\rho(X)=\rho(S)+1$,
\item $\dim S < \dim X$.
\endroster
\enddefinition

 Note that a Mori fibre space of relative dimension one is automatically a conic bundle.

\definition{Definition 2.5} Let $\map X.S.$ and $\map X'.S'.$ be Mori fibre spaces.  A
birational map $f:\rmap X.X'.$ is square if it induces a commutative diagram 
$$
\diagram
X     & \rDashto^f     & X'   \\
\dTo  &            & \dTo \\
 S    & \rDashto^g & S'   \\
\enddiagram
$$
where $g$ is a birational map and $f$ induces an isomorphism of the generic fibres.  
\enddefinition

 We recall a result of Corti [6,Theorem 4.2], which generalises a deep result of Sarkisov [19], in the following theorem.

\proclaim{Theorem 2.6} Let $\map X.S.$ be a standard conic bundle.
Let ${\Delta}\subset S$ denote the discriminant curve of the conic bundle, and assume that
$4K_S+{\Delta}$ is quasi-effective.  If $\map X'.S'.$ is another Mori fibre space, then
every birational map $\varphi:\rmap X.X'.$ is square.
\endproclaim 

 Recall that a $\Bbb Q$-divisor is called quasi-effective if it is a limit of effective
divisors in $N^{1}(S)$.

The following Lemma is a minor modification of Iskovskikh [8], Lemma 3 and Lemma 4.  The
only difference is that we have a birational map between $S$ and $S'$ instead of a
birational morphism.

\proclaim{Lemma 2.7} Let $\pi:\map X.S.$  be a standard conic bundle over smooth rational 
surface $S$. Let $S'$ be a rational surface and $\pi':\map X'.S'.$ a Mori fibre space with
one dimensional fibres which is square birational to $\pi:\map X.S.$, that is, there is a
commutative diagram
$$
\diagram
X'          & \rDashto^f & X        \\
\dTo^{\pi'} &            & \dTo_\pi \\
S'          & \rDashto^g & S        \\
\enddiagram
$$
where $f$ and $g$ are birational.  
\roster
\item If $\pi':\map X'.S'.$ is also a standard conic bundle, then the arithmetic genera of the
discriminant curves of these two standard conic bundles are equal.
\item The arithmetic genus of the discriminant curve of $\pi':\map X'.S'.$ is at least the 
arithmetic genus of discriminant curve of the standard conic bundle $\pi:\map X.S.$.
\endroster
\endproclaim
\demo{Proof} By elimination of indeterminacy, we know that there is a surface $\bar S$ and a 
commutative diagram 
$$
\diagram
              &             & \bar{S}    &                 &     \\
              &\ldTo^{\psi} &            & \rdTo^{\varphi} &     \\
S'            &             & \rDashto^g &                 &  S  \\
\enddiagram
$$
where the morphisms $\psi$ and $\varphi$ are compositions of smooth blow ups.  Applying
Sarkisov [19] Proposition 2.4 to $\varphi:\map \bar{S}. S.$, we obtain a standard conic
bundle $\bar{\pi}:\map \bar{X}. \bar{S}.$ such that the diagram

$$
\diagram
\bar{X}          & \rDashto           & X          \\
\dTo^{\bar{\pi}} &                & \dTo_{\pi} \\
\bar{S}          & \rTo^{\varphi} & S          \\
\enddiagram
$$
is commutative.

By Iskovskikh [8] Lemma 4, we know that the arithmetic genera of the discriminant curves
of these two standard conic bundles $\bar{\pi}:\map \bar{X}. \bar{S}.$ and $\pi:\map
X. S.$ are equal.  Now the diagram

$$
\diagram
\bar{X}          & \rDashto    & X'          \\
\dTo^{\bar{\pi}} &             & \dTo_{\pi'} \\
\bar{S}          & \rTo^{\psi} & S'          \\
\enddiagram
\tag \dag
$$
is also commutative (where the birational map $\rmap\bar{X}.X'.$ is the composition of
$\map \bar{X}.X.$ and the inverse of $f$).  The same argument shows that the arithmetic
genera of the discriminant curves of the standard conic bundles $\bar{\pi}:\map
\bar{X}. \bar{S}.$ and $\pi':\map X'.S'.$ are equal.  This completes the proof of part (1).

 Now we follow the proof of Iskovskikh [8] Lemma 3 to prove that the arithmetic genus of
the discriminant curve $C' \subset S'$ of the conic bundle $\pi':\map X'.S'.$ is at least
the arithmetic genus of the discriminant curve $\bar C \subset \bar S$ of the standard
conic bundle $\bar{\pi}:\map \bar X.\bar S.$.

Let $\Cal I_{\bar C}$ be the ideal sheaf of $\bar C$ in $\bar S$.  Then we
have an exact sequence of sheaves
$$
\ses \Cal I_{\bar C}.{\ring \bar S.}.{\ring \bar C.}..
$$

 If we push this exact sequence down to $S'$, then we get a long exact sequence, part of
which will be
$$
\es {0=R^1 \psi_*(\ring {\bar S}.)}.{R^1 \psi_*\ring \bar C.}.R^2 \psi_* \Cal I_{\bar C} = 0..
$$
We have zero on the right-hand side because the fibres of $\psi$ are at most one
dimensional and also on the left-hand side because the singularities of $\bar S$ and $S'$
are log terminal, whence rational.  The morphism $\psi$ induces a contraction morphism
$\psi |_{\bar C}: \map \bar C.C'.$, which is an isomorphism over a dense open subset of
$C'$, as $(\dag)$ is a commutative diagram.  Thus the Leray spectral sequence for $\psi
|_{\bar C}$ degenerates at the $E_2$-level and in particular $h^i (\Cal O_{\bar C})= h^i
(\psi_* \Cal O_{\bar C})$.

 The exact sequence 
$$
\ses {\ring C'.}.\psi_* {\ring \bar C.}.\Cal K.,
$$
where $\Cal K$ is defined by exactness, yields an exact cohomology sequence
$$
\ses H^0 (\Cal K).H^1  (\Cal O_ {C'}).H^1 (\psi_* \Cal O_ {\bar{C}})..
$$
This implies that $h^1(\Cal O_{\bar{C}})=h^1(\psi_*\Cal O_{\bar{C}}) \leq h^1(\Cal
O_{C'})$.  Combining with the first part of this Lemma gives the proof of part (2). \qed\enddemo

 We now want to focus on a sequence of conic bundles over $\pr 2.$, first introduced by 
Sarkisov [19].  

\proclaim{Lemma 2.9} Let $\pi:\map X.{\pr 2.}.$ be the projectivisation of the 
vector bundle $E=\ring {\pr 2.}.(k) \oplus \ring {\pr 2.}.\oplus \ring {\pr 2.}.$ over
$\pr 2.$.  Let $L$ be the tautological line bundle on $\proj E.$ so that $\pi_*L=E$ and
let $F$ be the pullback of the generator of $\Pic(\pr 2.)$.

 Then $\pl L.2.\otimes F$ is ample and base-point free and a generic divisor $V_k$ from
the linear system $|\pl L.2.\otimes F|$ on $X$ is a standard conic bundle over $\Bbb P^2$
whose discriminant curve is smooth of degree $2k+3$.
\endproclaim
\demo{Proof} The projectivisation $X$ of the 
vector bundle $E=\ring {\pr 2.}.(k) \oplus \ring {\pr 2.}.\oplus \ring {\pr 2.}.$ over
$\pr 2.$ has Picard number two and its cone of curves is generated by the class of 
a line contained in a fibre and the class of a line contained in the unique section with
normal bundle $\ring {\pr 2.}.(-k)\oplus \ring {\pr 2.}.(-k)$.  The line bundle $\pl
L.2.\otimes F$ is positive on both of these rays and so $\pl L.2.\otimes F$ is certainly
ample.  On the other hand $X$ is toric (indeed it is the projectivisation of a direct sum
of line bundles over a toric variety) and so $\pl L.2.\otimes F$ is automatically base-point free (e.g. see Fulton [7] page 70).

 Let $V_k$ be defined by $v_k = 0$ for some $v_k\in H^0(\pl L.2.\otimes F)$. Since
$$
\align
\pi_*(\pl L.2.\otimes F) &= \Sym ^2E \otimes \ring {\pr 2.}.(1)\\
                   &=\ring {\pr 2.}.(2k+1)\oplus\ring {\pr 2.}.(1)\oplus\ring {\pr 2.}.(1)\oplus\ring {\pr 2.}.(k+1)\oplus\ring {\pr 2.}.(1)\oplus\ring{\pr 2.}.(k+1),
\endalign
$$
we may write
$$
\pi_*(v_k) = (a_{11}, a_{22}, a_{33}, 2 a_{12}, 2 a_{23}, 2 a_{31})
$$
where $a_{11}\in H^0(\ring {\pr 2.}. (2k+1))$, $a_{22}\in H^0(\ring {\pr 2.}. (1))$,
$a_{33} \in H^0(\ring {\pr 2.}.  (1))$, $a_{12}\in H^0(\ring {\pr 2.}. (k+1))$, $a_{23}\in
H^0(\ring {\pr 2.}. (1))$ and $a_{31}\in H^0(\ring {\pr 2.}.(k+1))$.  It is easy to see
that the discriminant curve $\Delta(V_k)$ is the vanishing locus of
$$
D(\{a_{ij}\}) = \det 
\pmatrix a_{11} & a_{12} & a_{13}\cr
a_{21} & a_{22} & a_{23}\cr
a_{31} & a_{32} & a_{33}
\endpmatrix \in H^0(\ring {\pr 2.}. (2k+3))           \tag 1
$$
where we set $a_{ij}$ = $a_{ji}$.  In particular $\deg \Delta(V_k) = 2k+3$.

 It remains to show that for a general choice of $a_{ij}$, the curve $D(\{a_{ij}\})=0$
given in (1) is smooth. If $a_{12}=a_{23}=a_{31}=0$, the corresponding curve
$D(\{a_{ij}\})=0$ is a union $C_1\cup C_2\cup C_3$, where $C_1 =
\{a_{11} = 0\}$ is a general curve of degree $2k+1$ and $C_i=a_{ii}=0$, $i=2,3$ are two
lines in general position.  So the singularities of $D(\{a_{ij}\})=0$ are $C_i\cap C_j$.
We will show that these singularities are smoothed out as $v_k$ deforms.

Consider the family of curves given by
$$ 
\det 
\pmatrix a_{11} & {t a_{12}} & 0\cr
{t a_{12}} & a_{22} & 0\cr
0 & 0 & a_{33} \cr
\endpmatrix = (a_{11} a_{22} - t^2 a_{12}^2) a_{33} = 0.
$$ 

When $t = 0$, we have $C_1\cup C_2 \cup C_3$.  When $t\ne 0$, we see that the
singularities $C_1\cap C_2$ are smoothed out as long as we pick $a_{12}\in H^0(\ring {\pr
2.}. (k+1))$ not vanishing at $C_1\cap C_2$. Similarly, one can show that the
singularities $C_2\cap C_3$ and $C_3\cap C_1$ are smoothed out as $v_k$ deforms.
\qed\enddemo

\remark {Remark} The degree of the discriminant curve of $V_k$ was already known.  
It is well known that the discriminant curve of a standard conic bundle is nodal.  However
the fact that $\Delta(V_k)$ is smooth, when $V_k$ is general, is a new result due to Xi Chen. James
M$^c$Kernan even wonders if the following is true:

Let $X$ be the projectivisation of a rank 3 vector bundle over a smooth surface
$S$. Suppose that $|\pl L.2.\otimes \pl F.k.|$ is very ample, where $L$ restricts to
$\ring {\pr 2.}.(1)$ on a fibre.  Then for generic choice of $V$ belonging to the linear
system $|\pl L.2.\otimes \pl F.k.|$ the corresponding discriminant curve is smooth.
\endremark

\proclaim{Lemma 2.9} Let $k$, $m\geq 5$. If $V_k$ is birational to $V_m$, then $k=m$.
\endproclaim 
\demo{Proof} Assume $V_k$ is birational to $V_m$ and let $\Phi$ denote the
birational map.  As $k\geq 5$ Theorem 2.6 applies, therefore $\Phi:\rmap V_k. V_m.$ is
square, that is $\Phi$ induces a commutative diagram
$$
\diagram
V_k     & \rDashto^\Phi     & V_m   \\
\dTo  &            & \dTo \\
\pr 2.   & \rDashto & \pr 2..   \\
\enddiagram
$$
By Lemma 2.7 (1) the arithmetic genera of the discriminant curves of $V_k/\pr 2.$ and
$V_m/\pr 2.$ are equal.  Since the discriminant curves of $V_k$ and $V_m$ are smooth plane
curves of degree $2k+3$ and $2m+3$, respectively, we have that $k=m$. \qed\enddemo

\proclaim{Lemma 2.10} The standard conic bundle $V_k$ is birational to a $\Bbb Q$-Fano threefold $V_k'$ with Picard
number one.  
\endproclaim 
\demo{Proof} Consider the $k$-Uple Veronese embedding of $\Bbb P^2$ into $\Bbb P^N$, where 
$N=\frac{(k+2)(k+1)}{2}-1$.  Pick a line $l$ skew to $\pr N.$ in $\pr {N+2}.$.  Let $X'_k$
be the linear join of the image $S_k$ of $\pr 2.$ under this embedding and the line
$l$. Blowing up $X'_k$ along $l$ we obtain a fourfold $X_k$ which is a $\pr 2.$-bundle
over $\pr 2.$.  Now $X_k$ is isomorphic to $\Bbb P_{\pr 2.}(G)$, where $G$ is a rank $3$
vector bundle over $\pr 2.$. 

 The normal bundle $S_k$ in $X_k$ is $\Cal O_{\pr 2.} (k) \oplus \Cal O_{\pr 2.} (k) $.
On the other hand $S_k$ is a section of the natural map $\map X_k.{\pr 2.}.$ and so $S_k$
corresponds to a quotient $\map G.Q={\ring {X_k}.}(a).$ of $G$.  Let $K$ be the kernel.
Then the normal bundle of $S_k$ in $X_k$ is canonically isomorphic to $\Hom (K,
Q)=K^*\otimes Q$.  It follows that $K$ is split so that $K=\ring {\pr 2.}.(b)\oplus\ring
{\pr 2.}.(c)$ for some $b$ and $c$.  We have $a-b=a-c=k$.  Thus $b=c$.  Tensoring by a
line bundle, we may as well assume that $b=c=0$.  In this case $a=k$.  Thus
$X_k$ is indeed isomorphic to $\Bbb P_{\pr 2.}(E)$, where $E$ is the vector bundle
appearing in (2.8).
 
 Let $W_k$ be the image of $V_k$ inside $X_k'$ and let $f_k\:\map V_k.V_k'.$ be the
induced contraction, so that $V_k'$ is the normalisation of $W_k$.  Now $f_k$ is
birational; let $E$ be the exceptional locus (we use the same notation as we used for the
underlying vector bundle; hopefully this will not cause confusion).  Then $E$ is a divisor,
which intersects the general fibre $F$ of $\map V_k.{\pr 2.}.$ in two points.

 The Picard number of $X_k$ is clearly equal to two.  Hence $V_k$ also has Picard number
two, as $V_k$ is ample in $X_k$.  Thus $V_k'$ has Picard number one.  In particular $f_k$
has relative Picard number one.
 
As $X_k'$ has $\Bbb Q$-factorial singularities it follows that $K_{V_k'}$ is $\Bbb
Q$-Cartier.  Let $a$ be the log discrepancy of $E$, so that
$$
K_{V_k} + E = \pi^* K_{V_k'}+ a E.
$$ 
To compute $a$, pick a line $L$ inside $E$ and intersect both sides of this
equation with $L$,
$$ 
-3=K_E\cdot L= (K_{V_k} + E)\cdot L =(\pi^* K_{V_k'}+ a E)\cdot L=aE\cdot L=-ak.
$$ 
Thus $a=\frac3k$ and $V_k$ is log terminal, of log discrepancy $\frac 3k$.  Moreover 
it follows that 
$$
K_{V_k}+(1-a')E
$$ 
is $f_k$-negative, for any $a'>a$.  As $f_k$ has relative Picard number one, it follows
that $f_k$ is a step of the $K_{V_k}+(1-a')E$-MMP and in particular $V_k'$ is $\Bbb
Q$-factorial.

 Finally we check that $-K_{V_k'}$ is ample.  As $V_k'$ has Picard number one, it is
enough to check that $K_{V_k'}\cdot C<0$ for any curve $C$ in $V_k'$.  Let $F$ be a
general fibre of $\map V_k.{\pr 2.}.$ and let $C$ be the image of $F$ in $V_k'$.  Then
$K_{V_k}\cdot F=-2$ and $E\cdot F=2$.  Thus $K_{V_k}+(1-a)E$ is negative on $F$.  It
follows by push-forward that $K_{V_k'}\cdot C<0$, as required.  Thus $V_k'$ is indeed a $\Bbb
Q$-Fano threefold of Picard number one. \qed\enddemo

\subhead {3. Proof of Theorem 1.2}\endsubhead

Let $k$, $m\geq 5$ from now on.  Assume that the family of $\Bbb{Q}$-Fano threefolds
with Picard number one is birationally bounded.  Then there exists a parameter space $B$
of finite type and a morphism $\pi\:\map \Cal X.B.$ such that any $\Bbb Q$-Fano threefold is
birational to at least one fibre of $\pi$.  

 Let $S\subset B$ be the set of those points of $b\in B$ such that $X_b$ is birational to
$V_k$, for some $k$.  Now by Lemma 2.10, $V_k$ is birational to $V_k'$, a $\Bbb Q$-Fano
threefold, and moreover by Lemma 2.9, $V_k$ is not birational to $V_m$, for $m\neq k$.  It
follows that there are infinitely many fibers over $S$ which are not birational to each other.

 In the course of the proof of Theorem 1.2, we are going to repeatedly modify $B$
whilst still preserving the property that there are infinitely many fibers over $S$ which are not birational to each other.  For example we are clearly
free to replace $B$ by the closure of the points corresponding to $V_{k'}$ for infinitely many $k$.  The fact that there are infinitely many fibers over $S$ which are not birational to each other implies that the set of $k$ such 
that $V_k$ is birational to $X_b$ is infinite.

Let $Z$ denote the closure of $S$.  Then there is an irreducible component of $Z$ which
contains an infinite subset of $S$.  Renaming this component $B$, we may assume that the
set $S$ is dense in $B$ and that $B$ is irreducible.  Note that $B$ has positive
dimension, as $S$ is infinite.

 For the generic point $\eta \in B$, pick a resolution of $\Cal X_\eta$.  Replacing $B$ by
an open neighborhood of $\eta$ we may assume that $\Cal X / B$ is a flat family of smooth
projective three dimensional algebraic varieties over $B$.  We recall a deep result due to
Koll\'ar and Mori, [13], which states that one can run the MMP in families.

\proclaim{Theorem 3.1} Let $B$ be a connected normal quasi-projective variety and let 
$X/B$ be a flat, projective family of threefolds such that every fibre has only $\Bbb
Q$-factorial terminal singularities. 

 Then there is a finite, \'etale and Galois base change $p:\map B'.B.$, a flat projective
family $Y/B'$, and a rational map $f':\rmap X'/B'.Y/B'.$ such that on each fibre $f'$
induces a birational map, each fibre of $Y/B'$ has only $\Bbb Q$-factorial terminal
singularities and either 
\roster 
\item"(a)" $K_{Y/B'}$ is relatively nef, or  
\item"(b)" There is a morphism $\map Y.Z.$ over $B'$, where $-K_Y$ is relatively ample, and 
$\map Y_{b'}.Z_{b'}.$ is extremal, for every $b'\in B'$.
\endroster 
 Here $X'$ is the family over $B'$ obtained by base change.
\endproclaim
\demo{Proof} This is [13, Proposition 12.4.2].  \qed\enddemo

 Applying Theorem 3.1 to our situation, possibly passing to an \'etale cover, we may
assume that either (a) or (b) of Theorem 3.1 holds.  (a) is clearly impossible, since
infinitely many fibres are uniruled, and if $X_b$ is uniruled, then $K_{X_b}$ is certainly
not nef.

 Thus we may assume that there is a morphism $\pi:\map \Cal X.\Cal Z.$ over $B$, such that
$\map X_b.Z_b.$ is a Mori fibre space, for all $b\in B$.  Suppose that $b\in S$. As $k\geq
5$ Theorem 2.6 applies and the resulting birational map $\rmap X_b.V_k.$ is square.  In
particular for every $b\in S$, $Z_b$ is a surface.  As $S$ is dense in $B$, $\map \Cal
Z.B.$ is therefore a family of surfaces.  Equivalently, $\map \Cal X.\Cal Z.$ is a family
of conic bundles.

 Let $\Delta\subset \Cal Z$ be the locus where $\pi$ is not smooth.  Then for every point
$b\in S$, $\Delta_b$ is a curve in $Z_b$.  Thus, possibly passing to an open subset of
$B$, we may assume that the discriminant locus of $\map X_b.Z_b.$ is non-empty and equal
to $\Delta_b$.  Possibly passing to an open subset of $B$, we may assume that $\Delta$ is
flat over $B$ so that the arithmetic genera of the discriminant curves are constant over
$B$.

On the other hand, the arithmetic genus of the discriminant curve of the standard conic
bundle $V_k$ over $\pr 2.$ is $(k+1)(2k+1)$.  Suppose that $b\in S$.  According to Lemma
2.7 (2), the arithmetic genus of $\Delta_b$ is at least $(k+1)(2k+1)$.  As there are
infinitely many fibers over $S$ which are not birational to each other, the set of such
$k$ is infinite, and this contradicts the fact that the arithmetic genera are constant.

\centerline{\bf References}

\bigskip

\item{[1]} V. Alexeev. Boundedness and $K^2$ for log surfaces. Internat. J. Math. {\bf5} (1994), no.6, 779-810.

\item{[2]} A. Borisov. Boundedness theorem for Fano log-threefolds. J. Alg. Geom. {\bf5}(1)
(1996), 119-133.

\item{[3]} A. Borisov. Boundedness of Fano threefolds with log-terminal singularities of given index. J. Math. Sci. Univ. Tokyo {\bf 8} (2001), no.2, 329-342.

\item{[4]} A. Borisov and L. Borisov. Singular toric Fano varieties. Acad. Sci. USSR Sb. 
Math. {\bf75} (1993), no.1, 277-283.

\item{[5]} F. Campana. Une version g\'eom\'etrique g\'en\'eralis\'ee du th\'eor\`eme de produit 
de Nadel. Bull. Soc. Math. France {\bf 119} (1991), 479-493.

\item{[6]} A. Corti and M. Reid. Singularities of linear systems and $3$-fold birational geometry. 
In Explicit birational geometry of 3-folds. Cambridge University Press, 2000.

\item{[7]} W. Fulton. Introduction to toric varieties. Annals of Mathematics Studies, Vol. 131, Princeton University Press, Princeton, NJ, 1993, xii+157.

\item{[8]} V.A. Iskovskikh.  On a rationality criterion for conic bundles. Mat. Sb. {\bf187}(7) (1996), 
75-92.

\item{[9]} Y. Kawamata. Boundedness of $\Bbb Q$-Fano threefolds. Contemp. Math., Part 3, Amer. Math. Soc.,
{\bf131}, 1992.

\item{[10]} S. Keel and J. M$^c$Kernan. Rational curves on quasi-projective surfaces. Mem. Amer. Math. Soc.
 {\bf140} (1999), no.669.

\item{[11]} J. Koll\'ar, Y. Miyaoka and S. Mori. Rationally connectedness and boundedness of Fano manifolds. J. Diff. Geom., {\bf36} (1992), 765-779.

\item{[12]} J. Koll\'ar, Y. Miyaoka, S. Mori and H. Takagi. Boundedness of canonical $\Bbb Q$-Fano 3-folds. Proc. Japan Acad. Ser. A, Math. Sci., {\bf76}(5) (2000), 73-77.

\item{[13]} J. Koll\'ar and S. Mori. Classification of three-dimensional flips. J. Amer. Math. Soc. {\bf5} (1992), no. 3, 533-703.

\item{[14]} J. Koll\'ar and S. Mori. Birational geometry of algebraic varieties. Cambridge Tracts in Mathematics, 134, Cambridge University Press, Cambridge, 1998.

\item{[15]} J. M$^c$Kernan. Boundedness of log terminal Fano pairs of bounded index. 

arXiv:math.AG/0205214.

\item{[16]} A.M. Nadel. The boundedness of degree of Fano varieties with Picard number one. J. Amer. Math. Soc. {\bf4} (1991), 681-692.

\item{[17]} V. V. Nikulin. Del Pezzo surfaces with log terminal singularities. Mat. Sb. {\bf180} (1989), no.2, 226-243.

\item{[18]} Z. Ran and H. Clemens. A new method in Fano geometry. Internat. Math. Res. Notices, {\bf10}: 527-549, 2000.

\item{[19]} V. G. Sarkisov.  Birational automorphisms of conic bundles. Izv. Akad. Nauk SSSR Ser. Mat., {\bf44}(4) (1980), 918-945. (English translation: Math. USSR-Izv. {\bf17} (1981), no.4, 177-202.)

\end